\documentclass{amsart}
\usepackage{amssymb,latexsym}
\theoremstyle{plain}
\newtheorem{theorem}{Theorem}

\newtheorem{proposition}{Proposition}

\theoremstyle{definition}

\newtheorem{remark}{Remark}
\newtheorem{question}{Question}
\date{}

\begin{document}

\title[Hirzebruch surface]
{Vector bundles on Hirzebruch surfaces whose twists by a non-ample line 
bundle have
natural cohomology}
\author{E. Ballico}
\address{Dept. of Mathematics\\
 University of Trento\\
38050 Povo (TN), Italy}
\email{ballico@science.unitn.it}
\thanks{The author was partially supported by MIUR and GNSAGA of INdAM 
(Italy).}
\author{F. Malaspina}
\address{Dip. Matematica, Universit\`{a} di Torino\\
via Carlo Alberto 10, 10123 Torino, Italy}
\email{francesco.malaspina@unito.it}

\begin{abstract}
Here we study vector bundles $E$ on the Hirzebruch surface $F_e$ such 
that their twists by a spanned, but not
ample, line bundle $M = \mathcal {O}_{F_e}(h+ef)$ have natural 
cohomology, i.e. $h^0(F_e,E(tM)) >0$
implies $h^1(F_e,E(tM)) = 0$. 
\end{abstract}

\maketitle

\section{Introduction}\label{S1}

Let $F_e$, $e>0$, denote the Hirzebruch surface with a section with 
self-intersection $-e$. For any $L\in \mbox{Pic}(F_e)$
and any vector bundle $E$ on $F_e$ we will say that $E$ has property 
\pounds \pounds \    (resp. \pounds\ ) with respect
to $L$ if $h^1(F_e,E\otimes L^{\otimes m})=0$ for all $m \in \mathbb 
{Z}$ (resp. for all $m\in \mathbb {Z}$
such that $h^0(F_e,E\otimes L^{\otimes m}) \ne 0$). We think that 
property \pounds\  is nicer for
reasonable $L$. We take as a basis of $\mbox{Pic}(F_e) \cong \mathbb 
{Z}^2$
a fiber $f$ of the ruling $\pi : F_e \to {\bf {P}}^1$ and the section 
$h$ of $\pi$ with negative self-intersection.
Thus $h^2 = -e$, $h\cdot f = 1$ and $f^2 = 0$. We have $\omega _{F_e} 
\cong \mathcal {O}_{F_e}(-2h-(e+2)f)$.
$\mathcal {O}_{F_e}(\alpha h+\beta f)$ is spanned (resp. ample) if and 
only if $\alpha \ge 0$
and $\beta \ge \alpha e$ (resp. $\alpha > 0$ and $\beta > e\alpha$). 
The
Leray spectral sequence of $\pi$ and Serre duality give that 
$h^1(F_e,\mathcal {O}_{F_e}(\gamma h+\delta f))=0$
if and only if either $\gamma \ge 0$ and $\delta \ge e\gamma -1$ or 
$\gamma =-1$
or $\gamma \le -2$ and $-\delta - e-2 \ge e(-\gamma -2)-1$ (i.e. 
$\delta \le e\gamma +e-1$).
We consider as the test line bundle the spanned, but not ample, line 
bundle
$M:= \mathcal {O}_{F_e}(h+ef)$. Notice that the linear system $\vert 
M^{\otimes 2}\vert$ contains the sum of the effective divisor
$h$ and the ample divisor $h+2ef$. Thus for every vector bundle $E$ on 
$F_e$ there is an integer $m_0(E)$
such that $h^0(F_e,E\otimes M^{\otimes m}) \ne 0$ for all $m \ge 
m_0(E)$. We will see that property \pounds \ \pounds\ 
is too strong and not interesting (see Remarks \ref{a0} and \ref{a1}). 
We stress the property \pounds\  with respect
to $M$ is quite different from similar looking properties (e.g. natural 
cohomology) with respect to an ample line bundle.\
(see Remarks \ref{a1} and \ref{a2} for the rank $1$ case).
Obviously, properties \pounds\  and \pounds \ \pounds\ 
may be stated for arbitrary projective varieties. In dimension $n \ge 
3$, one need to choose between vanishing
of $h^1$ or vanishing of all $h^i$, $1 \le i \le n-1$. We considered 
here the example $(F_e,M)$, because
it is geometrically significant. Indeed, let $\phi _M$ denote the 
morphism associated to the
base point free linear system $\vert M\vert$. If $e=1$ the morphism 
$\phi _M$ is the blowing up $F_1 \to {\bf {P}}^2$.
If $e\ge 2$, then $\phi _M : F_e \to {\bf {P}}^{e+1}$ contracts $h$ and 
its image is a cone over the rational normal
curve of ${\bf {P}}^e$. Moreover, for any spanned and non-trivial line 
bundle $L$ on $F_e$ there is an effective divisor
$D$ such that $L \cong M(D)$. For any spanned, but not ample line 
bundle $A$ on $F_e$ there is an integer
$c \ge 0$ such that $A \cong M^{\otimes c}$. 
We prove the following results.

\begin{theorem}\label{o4}
Fix integers $e \ge 1$, $r \ge 1$, $u, v$ such that $v \le e(u-r+1)-2$. 
Then there is no rank $r$ vector bundle
$E$ on $F_e$ with property \pounds\  with respect to $M$ such that 
$c_1(E) = \mathcal {O}_{F_e}(uh+vf)$.
\end{theorem}

\begin{theorem}\label{o2}
Fix integers $e, m, u, v$ such that $e\ge 1$, and $v \ge e(u-1)-1$ and 
$m \ge 0$. Set $\widetilde{a}:=
\sum _{i=0}^{u+2m-2} v+2m-1-ie$
and $\widetilde{b}:= \sum _{i=0}^{u+2m-1} v+2m -ie$. Fix any integer 
$s$ such that $\widetilde{a}
\le s \le \widetilde{b}$. Then there exists a rank $2$ vector bundle
$E$ on $F_e$ with property \pounds\  with respect to $M$ such that 
$c_1(E) = \mathcal {O}_{F_e}(uh+vf)$
and $c_2(E):=s -e(u+m-1) + (1-m)(v+em)$. Set $R:= \mathcal 
{O}_{F_e}(h+(e+1)f)$
and assume $m=0$, $u \ge 3$, and $v < 2eu$. Then
we may find $E$ as above which is $R$-stable in the sense of 
Mumford-Takemoto and (under the additional condition
$v \le 2eu -3$) such that $N\cdot M < c_1(E)\cdot M/2$
for all rank $1$ subsheaves $N$ of $E$.
\end{theorem}

The case `` $r=1$ '' of Theorem \ref{o4} is obviosly true (use the 
cohomology of line bundles on $F_e$, i.e. Remark
\ref{a1} below). In this case the converse is true, i.e. $\mathcal 
{O}_{F_e}(uh+vf)$ has property \pounds\ with
respect to $M$ if and only if $v \ge eu-1$ (Remark \ref{a1}). We were 
surprised
thar for $r \ge 2$ there is no way to overcome this $c_1$-obstruction.

The assumptions of the last part of Theorem \ref{o2} may be relaxed and 
instead of $R$ we may take
an arbitrary ample divisor $H$. An interesting offshot of our proof of 
Theorem \ref{o2} is that
our examples are given by an extension (\ref{eqo9}) and all locally 
free sheaves fitting in (\ref{eqo9})
have property \pounds\  with respect to $M$ and (under the additional 
conditions listed in Theorem \ref{o2})
are $R$-stable and $N\cdot M < c_1(E)\cdot M/2$
for all rank $1$ subsheaves $N$ of $E$.

In the case of direct sums of line bundles we will prove the following 
result.

\begin{proposition}\label{o3}
Fix integers $e \ge 1$, $r \ge 2$ and $L_i\in \mbox{Pic}(F_e)$, $1 \le 
i \le r$, say $L_i \cong
\mathcal {O}_{F_e}(u_ih+v_if)$. Set $E:= L_1\oplus \cdots \oplus L_r$. 
Up to a permutation of the factors
of $E$ we may assume $u_1 \ge \cdots \ge u_r$ and that if $u_i = u_j$ 
for some $i<j$, then $v_i \ge v_j$.
Set $m:= -u_1$ if $v_1 \ge eu_1$ and $m:= -u_1+1$ if $v_1 = eu_1-1$. 
The vector bundle
$E$ has property \pounds\  with respect to $M$ if and only if $v_i \ge 
eu_i-1$ for all $i$, and
for each $i \in \{2,\dots ,r\}$
either $u_i-m \ge -1$ or $-1 \le v_i-eu_i \le e-1$.
\end{proposition}

We raise the following question.

\begin{question}\label{o3.1}
Assume $e=1$ or $e=2$. Is it possible to describe all invariants $r, 
c_1, c_2$ of vector bundles on $F_e$
with property \pounds \ with respect to $M$?
\end{question}

\section{The proofs}\label{S2}

For any sheaf $F$ we will often write $F(mM)$ instead of $F\otimes 
M^{\otimes m}$.

\begin{remark}\label{a0}
The line bundle $\mathcal {O}_{F_e}(ch+df)$ is ample if and only if 
$c>0$ and $d > ec$. Hence any ample line bundle
is spanned. Assume that $H:= \mathcal {O}_{F_e}(ch+df)$ is ample. The 
cohomology of line bundles on $F_e$ shows
that for every $t\in \mathbb {Z}$ the line bundle $H^{\otimes t}$ has 
property \pounds \ \pounds\  with respect to $H$.
Hence $\mathcal {O}_{F_e}$ has property \pounds \ \pounds\  with 
respect to any ample line bundle. Set $H':= \mathcal {O}_{F_e}(ch+(d+2)f)$.
Taking $m_t:= -tc$ we see that no $H^{\otimes t}$, $t > 0$, has 
property \pounds \ \pounds\  with respect to the ample
line bundle $H'$. Taking $m_t:= -tc$ we see that no $H^{\otimes t}$, $t 
> 0$, has property \pounds \ \pounds\  with respect to $M$.
\end{remark}

\begin{remark}\label{a1}
Here we  study properties \pounds\  and \pounds \ \pounds\  with 
respect to $M$ for line bundles on $F_e$. Fix
$L\in \mbox{Pic}(F_e)$, say $L \cong \mathcal {O}_{F_e}(u h + vf)$. 
First assume $v \ge eu$.
We have $h^0(F_e,L(xM)) >0$ if and only if $x \ge -u$. Since 
$h^1(F_e,\mathcal {O}_{F_e}(ch+df)) = 0$
if $c \ge 0$ and $d \ge ec$, $L$ has property \pounds\  with respect to 
$M$.
Now assume $v < eu$. We have $h^0(F_e,L(xM)) > 0$ if and only if $ex 
\ge -v$. Since $h^1(F_e,\mathcal {O}_{F_e}(ch+df)) = 0$
if $c \ge 0$ if and only if $d \ge ec-1$, we get that $L$ has property 
\pounds\  with respect to $M$
if and only if $v \ge eu-1$. Take $m:= -u$. If $v=eu$, then we saw in 
the
introduction that $L$ has property \pounds \ \pounds\  with respect to 
$M$ if and only if $e=1$. Notice that
$h^1(F_e,\mathcal {O}_{F_e}((u-x)h + (v-ex)f)) = h^1(F_e,\mathcal 
{O}_{F_e}((x-u-2)h+(ex-v-e-2)) >0$
when $x \ge -u-2$ if and only if $-eu -2e \le -v-e-1$, i.e. if and only 
if $v \le eu+e-1$. Notice that 
$h^1(F_e,\mathcal {O}_{F_e}((u+x)h + (v+ex)f)) =0$ for $x \ge -u $ if 
and only if $v \ge eu-1$. Notice
that $h^1(F_e,\mathcal {O}_{F_e}(-h + (v-eu-e)f)) = 0$ for every $v\in 
\mathbb {Z}$. Hence
$L$ has property \pounds \ \pounds\  with respect to $M$ if and only if 
$eu - 1 \le v \le eu+e-1$. 
\end{remark}

\begin{remark}\label{a2}
Here we look at property \pounds\  with respect to the ample line 
bundle
$R:= \mathcal {O}_{F_e}(h+(e+1)f)$ for line bundles on $F_e$. Fix
$L\in \mbox{Pic}(F_e)$, say $L \cong \mathcal {O}_{F_e}(u h + vf)$. We 
heve $h^0(F_e,L(xR)) = 0$
if and only if $x \ge -u$ and $x(e+1) \ge -v$. We immediately see that 
if $v \ge (e+1)u$, then $L$ has
property \pounds\  with respect to $M$. Now assume $v < (e+1)u$. Set 
$y:= \lceil -v/(e+1)\rceil$.
We have $h^0(F_e,L(xR)) >0$ if and only if $x \ge y$. Fix an integer $x 
\ge y$.
Since $u+x \ge u+y \ge 0$,
$h^1(F_e,\mathcal {O}_{F_e}(u+x)h+(v+(e+1)x)f)) >0$ if and only if 
$v+(e+1)x \le eu+ ex-2$. The strongest condition
is obtained when $x=y$. We get that $L$ has property $\alpha $ with 
respect to $R$
if and only if either $v \ge (e+1)u$ or $v+ey \ge eu-1$, where $y:= 
\lceil -v/(e+1)\rceil$. 
\end{remark}

\begin{remark}\label{a3}
$E_1\oplus E_2$ has property \pounds \ \pounds\  with respect to $L$ if 
and only if both $E_1$ and $E_2$ have
property \pounds \ \pounds\  with respect to $L$. If $E_1\oplus E_2$ 
has property \pounds\  with respect
to $L$, then the same is true for $E_1$ and $E_2$. Now we check that 
the converse is not true.
Both $\mathcal {O}_{F_e}$ and $\mathcal {O}_{F_e}(-2h+(-e+4)f)$ have 
property \pounds\  with respect to $M$
(Remark \ref{a1}). Since $h^1(F_e,\mathcal {O}_{F_e}(-2h+(-e+4)f)) = 
h^1(F_e,\mathcal {O}_{F_e}(-2f)) = 1$,
$\mathcal {O}_{F_e}\oplus \mathcal {O}_{F_e}(-2h+(-e+4)f)$ has not 
property \pounds\  with respect to $M$.
\end{remark}

\begin{remark}\label{a3.0}
The definition of property \pounds\  may be given for an arbitrary 
torsion free sheaf, but not much may be said
in the general case. Here
we look at the rank $1$ case, because we will need it in the proofs of 
Theorems \ref{o4} and \ref{o2}. Let $A$ be a rank $1$ torsion free 
sheaf on $F_e$. Hence $A \cong \mathcal {I}_Z(uh+vf)$
for some zero-dimensional scheme $Z$ and some integers $u,v$. Since $Z$ 
is zero-dimensional,
$h^1(F_e,\mathcal {O}_{F_e}((u+t)h+(v+et)f)) \le h^1(F_e,\mathcal 
{I}_Z((u+t)h+(v+et)f))$ for all $t\in \mathbb {Z}$.
Taking $t \gg 0$ we see that if $A$ has property \pounds\  with respect 
to $M$, then $v \ge eu-1$.
When $v \ge eu-1$, for a general $Z$ (in the following sense) $A$ has 
property \pounds\  with respect to $M$
for the following reason. Fix an integer $z>0$. Since $F_e$ is a smooth 
surface, the Hilbert scheme $\mbox{Hilb}^z(F_e)$ of all length $z$ 
zero-dimensional
subschemes of $F_e$ is irreducible and of dimension $2z$ (\cite{f}). 
Take a general $S\in \mbox{Hilb}^z(F_e)$, i.e.
take $z$ general points of $F_e$. Since $h^0(F_e,\mathcal {I}_S\otimes 
L) =
\max \{0,h^0(F_e,L)-z\}$ for every $L\in \mbox{Pic}(F_e)$, it is easy 
to check that if $v \ge eu-1$, then
$\mathcal {I}_S(uh+vf)$ has property \pounds\  with respect to $M$. Now 
take $v = eu-1$, any integer $z>0$ and
any zero-dimensional length $z$ subscheme $B$ of $h$. Twisting with 
$(-u+1)M$ we see that $\mathcal {I}_B(uh,(eu-1)f)$
has not property \pounds\  with respect to $M$. Now assume $v > eu$. 
Take a zero-dimensional length $z\ge 2$ scheme
$W$ of a fiber of $\pi$. Twisting with $-uM$. We see that $\mathcal 
{I}_W(uh+vf)$ has not property \pounds\ 
with respect to $M$. If $z \ge 3$ and $v = eu$, twisting with $(-u+1)M$ 
and using the same $W$ we get a sheaf without
property \pounds\  with respect to $M$.
\end{remark}

Property \pounds\  with respect to $M$ has the following open property.

\begin{proposition}\label{a4}
Let $\{E_t\}_{t\in T}$ be a flat family of vector bundles on $F_e$ 
parametrized by an integral variety
$T$. Assume the existence of $s\in T$ such that $E_s$ has property 
\pounds\  with respect to $M$.
Then there exists an open neighborhood $U$ of $s$ in $T$ such that 
$E_t$ has property \pounds\  with
respect to $M$ for all $t\in U$.
\end{proposition}

\begin{proof}
Let $m$ be the minimal integer such that $h^0(F,e,E_s(mM)) >0$. Thus 
$h^1(F_e,E_s(xM)) = 0$ for all $x \ge m$.
By semicontinuity there is an open neighborhood $V$ of $s$ in $T$ such 
that $h^0(F_e,E_t((m-1)M)) = 0$ for all $t\in V$.
By semicontinuity for every integer $x \ge m$ there is an open 
neighborhhod $V_x$ of $s$ in $T$
such that $h^1(F_e,E_t(xM))$ for all $t\in V_x$. Fix an irreducible 
$D\in \vert M\vert$. Hence
$D \cong {\bf {P}}^1$. Since $D^2 >0$, there
is an integer $a$ such that $h^1(D,E_s(aM)\vert D) = 0$. By 
semicontinuity there is an open
neighborhood $V$ of $s$ in $T$ such that $h^1(D,E_t(aM)\vert D) = 0$ 
for every $t\in V$. Since
$D^2>0$, $h^1(D,E_t(xM)\vert D) = 0$ for every $t\in V$ and every 
integer $x \ge a$. Fix an integer
$x \ge a$. From the exact sequence
\begin{equation*}
0 \to E_t((x-1)M) \to E_t(xM) \to E_t(xM)\vert D \to 0
\end{equation*}
we get that if $h^1(F_e,E_t((x-1)M)) = 0$, then $h^1(F_e,E_t(xM))$. 
Hence we may
take $U:= V\cap \bigcap _{x=m}^{\max \{a,m\}} V_x$.
\end{proof}

\vspace{0.3cm}

\qquad {\emph {Proof of Proposition \ref{o3}.}} If $E$ has property 
\pounds\  with respect to $M$, then each $L_i$
has property \pounds\  with respect to $M$ (Remark \ref{a3}) and hence 
$v_i \ge eu_i-1$ for all $i$. Now we assume
$v_i \ge eu_i-1$ for all $i$. Notice that $m$ is the minimal integer 
$t$ such that $h^0(F_e,E(t)) \ne 0$.
Since $L_1$ has property \pounds\  with respect to $M$, $E$ has 
property \pounds\  with respect to $E$
if and only if $h^1(F_e,L_i(tM)) =0$ for all $t \ge m$ and all 
$i=2,\cdots ,r$. If $u_i-m \ge -1$, then $h^1(F_e,L_i(tM)) =0$ for all $t \ge 
m$
because $v_i \ge eu_i-1$. Now assume $u_i-m \le -2$. We get 
$h^1(F_e,L_i(tM)) =0$ for any $t \ge m$ if and
only if $-1 \le v_i-eu_i \le e-1$.\qed

\vspace{0.3cm}

Here we discuss the set-up for the rank $2$ case. Consider an exact 
sequence
\begin{equation}\label{eqo1}
0 \to \mathcal {O}_{F_e}(D) \to E(mM) \to \mathcal {I}_Z(c_1+2mM-D) \to 
0
\end{equation}
in which $Z$ is a zero-dimensional scheme with length $s$
and either $D=0$ or $D = h$ or $D \in \vert zf\vert$ for some $1 \le z 
\le e$
or $e \ge 2$ and $D \in \vert h+wf\vert$ for some $1 \le w \le e-1$. We 
have $c_1(E(mM)) = c_1+2nm$ and $c_2(E(mM))
= s + D\cdot c_1 +2mM\cdot D -D^2$. Thus $c_1(E) = c_1$ and $c_2(E) = 
c_2$ by the choice of $s$ (\cite{h2}, Lemma 2.1).
Each $E$ fitting in (\ref{eqo1}) is torsion free. To have some locally 
free $E$ fitting in (\ref{eqo1})
a necessary condition is that $Z$ is a locally complete intersection. 
Notice that $h^1(F_e,\mathcal {O}_{F_e}(D)(-M))
= 0$ if $h$ is not a component of $D$. Hence a sufficient condition to 
have $h^0(F_e,E(mM)) >0$ and $h^0(F_e,E((m-1)M))=0$
is the equality
\begin{equation}\label{eqo2}
h^0(F_e, \mathcal {I}_Z(c_1+(2m-1)M-D))=0
\end{equation}
and (\ref{eqo2}) is a necessary condition if $h$ is not a component of 
$D$. Assume that $Z$ is a locally a complete intersection. The 
Cayley-Bacharach
condition associated to (\ref{eqo1}) is satisfied if 
\begin{equation}\label{eqo3}
h^0(F_e,\mathcal {I}_{Z'}(c_1+2mM-2D-2h -(e+2)f)) =0
\end{equation}
for every length $s-1$ closed subscheme of $s$ (\cite{c}). This 
condition is satisfied
if $h^0(F_e, \mathcal {I}_Z(c_1+(2m-1)M-D))=0$, $Z_{red}\cap h = 
\emptyset$
and no connected component of $Z$ is tangent to a fiber of the fiber of 
$\pi$, because the line bundle
$\omega _{F_e}^\ast (-M)
= \mathcal {O}_{F_e}(h+2f)$ is base point free outside $h$ and the 
morphism associated to
$\vert f\vert$ is the ruling; if $e =12$, then (\ref{eqo3}) is 
satisfied if (\ref{eqo2}) is satified, because $\mathcal {O}_{F_1}(h+2f)$
is very ample; if $e=2$ it is sufficient to assume $Z_{red}\cap h = 
\emptyset$, because the
morphism associated to $\mathcal {O}_{F_2}(h+2f)$ is an embedding 
outside $h$. 

\vspace{0.3cm}

\qquad {\emph {Proof of Theorem \ref{o4} for $r \le 2$.}} If $r=1$, 
then use Remark \ref{a1}.
Assume the existence of a rank two vector bundle $E$ with
property \pounds\  with respect to $M$ and $c_1(E) = \mathcal 
{O}_{F_e}(uh+vf)$. Let $m$ be the first integer
such that $h^0(F_e,E(mM)) >0$. We get an exact sequence (\ref{eqo1}) 
with $D\in \vert \mathcal {O}_{F_e}(xh+yf)$
with the convention $(x,y)=(0,0)$ if $D = \emptyset$. Hence either 
$(x,y) = (0,0)$ or $(x,y) = (1,0)$
or $x=0$ and $1 \le y \le e$ or $e \ge 2$, $x=1$, and $1 \le y \le 
e-1$. Since $h^2(F_e,M^{\otimes z}(D)) = 0$ for
all $z \ge 0$, (\ref{eqo1}) and property \pounds\  for $E$ imply 
$h^1(F_e,\mathcal {I}_Z((u+2m-x+z)h+(v+2me-y+ze)f)=0$
for all $z \ge 0$. As in Remark \ref{a3.0} we see that when $z \gg 0$ 
the last equality
implies $v-y \ge e(u-x)-1$. If $v \le e(u-1)-2$ the last inequality is 
not satisfied for any choice of the pair
$(x,y)$ in the previous list.\qed 

\vspace{0.3cm}

\quad {\emph {Proof of Theorem \ref{o2}.}} Fix a general $S \subset 
F_e$
such that $\sharp (S)=s$. Let $E$ be any torsion free sheaf fitting in 
the following exact sequence:
\begin{equation}\label{eqo9}
0 \to \mathcal {O}_{F_e}((1-m)h-emf) \to E \to \mathcal 
{I}_S((u+m-1)h+(v+em)f) \to 0
\end{equation}
We have $c_1(E) = \mathcal {O}_{F_e}(uh+vf)$ and $c_2(E) = s -e(u+m-1) 
+ (1-m)(v+em)$. By construction
$h^0(F_e,E(mM)) \ne 0$. We have $h^0(F_e,E((m-1)M)) = 0$. If 
$h^0(F_e,\mathcal {I}_S((u+2m-2)h+(v+2em-e)f)=0$.
Since $S$ is general, $h^0(F_e,\mathcal {I}_S((u+2m-2)h+(v+2em-e)f)=0$ 
if and only if
\begin{equation}\label{eqo10}
h^0(F_e,\mathcal {O}_{F_e}((u+2m-2)h+(v+2em-e)f)\le s
\end{equation}
Since $S$ is general, every subset 
of it is general. Hence to check the Cayley-Bacharach condition and 
hence show the local freeness of a general $E$
given by the extension (\ref{eqo10})
it is sufficient to prove check the following inequality:
\begin{equation}\label{eqo11}
h^0(F_e,\mathcal {O}_{F_e}((u+2m-5)h+(v+2em-2e-2)f))\le s-1
\end{equation}
This is true, because we assumed $s \ge \widetilde{a}$ and 
$\widetilde{a} >
h^0(F_e,\mathcal {O}_{F_e}((u+2m-5)h+(v+2em-2e-2)f))$. Hence a general
$E$ fitting in the extension (\ref{eqo10}) is locally free. Since 
$\mathcal {O}_{F_e}(3h+e+2)$ has a a subsheaf the very ample line bundle 
$\mathcal {O}_{F_e}(h+e+2)$, (\ref{eqo11})
is satisfied if (\ref{eqo10}) is satisfied. The generality of $S$ 
implies that $h^1(F_e,\mathcal {I}_S((u+m-1+t)h+(v+em+et)f)=
0$ if and only if $h^1(F_e,\mathcal {O}_{F_e}((u+m-1+t)h+(v+em+et)f)=0$ 
and $h^0(F_e,\mathcal {O}_{F_e}((u+m-1+t)h+(v+em+et)f)
\ge s$. Notice that $\widetilde{a} = h^0(F_e,\mathcal 
{O}_{F_e}((u+2m-2)h
+ (v+2me -e)f))$ and $\widetilde{b} = h^0(F_e,\mathcal 
{O}_{F_e}((u+2m-1)h
+ (v+2me )f))$. Since $h^1(F_e,(u+m-1+t)u + (v+me+te)f)) = 0$
for all $t \ge 0$, $\widetilde{a} \le s \le \widetilde{b}$ and $S$ is 
general, any sheaf $E$ in (\ref{eqo9}) has
property \pounds\  with respect to $M$. Since a general extension 
(\ref{eqo9}) has locally free middle term $E$,
the proof of the first part of Theorem \ref{o2} is over. Now assume 
$m=0$, $u \ge 3$,$v < 2eu$, and that $E$ is not $R$-stable, i.e. assume
the existence of $N\in \mbox{Pic}(F_e)$ such that $N\cdot R \ge 
c_1(E)\cdot R/2$ and an inclusion $j:
N \to E$; here to have $N$ locally free we use that $E$ is reflexive. 
Since $m=0$ and $u \ge 3$, $c_1(E)\cdot R > 
2(\mathcal {O}_{F_e}(h))\cdot R$. Hence $j$ induces a non-zero map
$N \to \mathcal {I}_S((u-1)h+vf)$. Any non-zero map $N \to \mathcal 
{O}_{F_e}((u-1)h + vf)$
is associated to a unique non-negative divisor $\Delta \in \vert
\mathcal {O}_{F_e}((u-1)h + v)f)\otimes N^\ast \vert$. Since $j$
factors through $\mathcal {I}_S((u-1)h+vf)$, $h^0(F_e,\mathcal 
{I}_S(\Delta ) >0$. We fixed $R$ and
the integers $m, u, v$. There are only finitely many possibilities for 
the line bundle
$\mathcal {O}_{F_e}(\Delta )$. Since $S$ is general, we
get $h^0(F_e,\mathcal {O}_{F_e}(\Delta )) \ge s$. Write $N = \mathcal 
{O}_{F_e}(\gamma h+\delta f)$
for some integers $\gamma ,\delta$. The inequality $N\cdot R \ge 
c_1(E)\cdot R/2$ is equivalent to the
inequality
\begin{equation}\label{eqo12}
2\gamma + 2\delta \ge u+v
\end{equation}
We have $\mathcal {O}_{F_e}(\Delta ) = \mathcal {O}_{F_e}((u-1 -\gamma 
)h+(v -\delta )f)$.
Since $h^0(F_e,\mathcal {O}_{F_e}(\Delta )) \ge s $ and $s \le 
\widetilde{b} = h^0(F_e,\mathcal {O}_{F_e}((u-1)h
+ (v )f))$, either $\gamma \le 0$ or $\delta \le 0$. Since $\Delta$ is 
effective, we also
have $\gamma \le u-1$ and $\delta \le v$. First assume $\delta \le 0$. 
Hence $\gamma \ge (u+v)/2$.
Since $ \gamma \le u-1$, we get $v \le u-2$. Since $v \ge eu -e$, we 
get a contradiction.
Now assume $\gamma \le 0$. We get $\delta \ge (u+v)/2$. Consider the 
exact sequence
\begin{equation}\label{eqo13}
0 \to N \to E \to \mbox{Coker}(j) \to 0
\end{equation}
Notice that $\mbox{Coker}(j)^{\ast \ast }\cong \mathcal 
{O}_{F_e}((u-\gamma )h+(v -\delta )f)$.
Since $\gamma \le 0$, $\delta \ge (u+v)/2$, and $v < 2eu$, we have $v- 
\delta \le e(u -\gamma )-2$.
In Remark \ref{a3.0} we checked that $h^1(F_e,\mbox{Coker}(j)(tM)) > 0$ 
for $t \gg 0$. Since $h^2(F_e,L(tM))
=0$ for $t \gg 0$ and any $L\in \mbox{Pic}(F_e)$, the exact sequence 
(\ref{eqo13}) gives
that $E$ has not property \pounds\  with respect to $M$, contradicting 
the already proved part
of Theorem \ref{o2}. If instead of $R$ we use $M$ for the intersection 
product, instead of (\ref{eqo12})
we only have the inequality $2 \delta \ge v$. Everything works in the 
same way with only minor
numerical modifications.\qed

\begin{remark}\label{o5}
There are at least $2$ well-known and related ways to obtain rank $r 
\ge 3$ vector bundles as extensions.
Instead of (\ref{eqo1}) we may take the exact sequence
\begin{equation}\label{eqo4}
0 \to \oplus _{i=1}^{r-1} \mathcal {O}_{F_e}(D_i-m_iM) \to \mathcal 
{I}_Z(uh+vf) \to 0
\end{equation}
In \cite{hl}, proof of Theorem 5.1.6, the following extension is used:
\begin{equation}\label{eqo5}
0 \to L_1 \to E \to \oplus _{i=2}^{r} \mathcal {I}_{Z_i}(u_ih+v_if) \to 
0
\end{equation}
The latter extension was behind the proof of Proposition \ref{o3}. Both 
extensions can give several examples
of vector bundles with or without property \pounds\  with respect to 
$M$. To prove Theorem \ref{o4} we will
use iterated extensions, i.e. increasing filtrations $E_i$, $1 \le i 
\le r$, of $E$ such that $E_1$
is a line bundle, $E_r = E$ and each $E_i/E_{i-1}$ is a rank $1$ 
torsion free sheaf. 
\end{remark}

\quad {\emph {Proof of Theorem \ref{o4} for $r \ge 3$.}} Assume the 
existence of a rank $r$ vector bundle $E$ with
property \pounds\  with respect to $M$ and $c_1(E) = \mathcal 
{O}_{F_e}(uh+vf)$. Let $m_1$ be the first integer
such that $h^0(F_e,E(m_1M)) >0$. Fix a general $\sigma \in 
H^0(F_e,E(m_1M))$. Since $h^0(F_e,E((m_1-1)M)) =0$,
$\sigma$ induces an exact sequence
\begin{equation}\label{eqo6}
0 \to \mathcal {O}_{F_e}(-m_1M+D_1) \to E \to G_1 \to 0
\end{equation}
with $F_1$ torsion free, $D_1$ of type $(x_1,y_1)$ and either 
$(x_1,y_1) = (0,0)$ or $(x_1,y_1) = (1,0)$
or $x_1=0$ and $1 \le y_1 \le e$ or $e \ge 2$, $x_1=1$, and $1 \le y_1 
\le e-1$. Notice
that $c_1(G_1) = \mathcal {O}_{F_e}((u+m_1-x_1)h + (v+em_1-y_1)f)$. Set 
$E_1:= \mathcal {O}_{F_e}(-m_1M+D_1)$.
Since $h^2(F_e,\mathcal {O}_{F_e}((t-m_1)D+D_1)) = 0$ for all $t \ge 
m_1$, property \pounds
for $m$ with respect to $M$ implies $h^1(F_e,F_1(tM)) = 0$ for all $t 
\ge m_1$. Let $m_2$ be the first
integer such that $m_2 \ge m_1$ and $h^0(F_e,F_1(m_2M)) >0$. A non-zero 
section of $H^0(F_e,G_1(m_2M))$
induces an exact sequence
\begin{equation}\label{eqo7}
0 \to \mathcal {I}_{Z_1}(-m_2M+D_2) \to G_1 \to G_2 \to 0
\end{equation}
with $\mathcal {I}_{Z_1}$ zero-dimensional, $G_2$ torsion free and 
$D_2$ an effective divisor of type $(x_2,y_2)$ and either $(x_2,y_2) = 
(0,0)$
or $(x_2,y_2) = (1,0)$
or $x_2=0$ and $1 \le y_2 \le e$ or $e \ge 2$, $x_2=1$, and $1 \le y_2 
\le e-1$. Here we cannot claim that $Z_1 = \emptyset$,
because $G_1$ is not assumed to be locally free. Notice
that $c_1(G_2) = \mathcal {O}_{F_e}((u+m_1+m_2-x_1-x_2)h + 
(v+em_1+em_2-y_1-y_2)f)$.
Since $Z_1$ is zero-dimensional, $h^2(F_e,\mathcal {I}_{Z_1}\otimes L)= 
h^2(F_e,L)$
for every $L\in \mbox{Pic}(F_e)$. Hence as in the first step we get 
$h^1(F_e,G_2(tM)) = 0$ for all $t \ge m_2$.
If $r=3$, we are done as im the proof of the case $r=2$. If $r \ge 4$, 
we iterate the last step $r-3$ times.\qed

\providecommand{\bysame}{\leavevmode\hbox to3em{\hrulefill}\thinspace}

\end{document}